# On the Diagonalization of Difference Calogero-Sutherland Systems

J. F. van Diejen

ABSTRACT. We discuss the simultaneous diagonalization of a family of commuting difference operators by Koornwinder's multivariable generalization of the Askey-Wilson polynomials. The operators constitute a complete set of quantum integrals for a difference version of the $n$-particle quantum Calogero-Sutherland system related to the root system $BC_n$.

## 1. Introduction

Some years ago, Ruijsenaars introduced [**19**] an integrable quantum $n$-particle model that can be viewed as a relativistic generalization of the well-known quantum Calogero-Sutherland system [**22, 23**]. The commuting quantum integrals for the model are given by difference operators of the form

(1.1)
$$\hat{H}_r \;=\; \sum_{\substack{J\subset\{1,\ldots,n\} \\ |J|=r}} \prod_{\substack{j\in J \\ k\not\in J}} v^{\frac{1}{2}}(x_j - x_k)\, e^{i\beta \sum_{j\in J} \partial_j} \prod_{\substack{j\in J \\ k\not\in J}} v^{\frac{1}{2}}(x_k - x_j), \qquad r = 1,\ldots,n,$$

where $\partial_j = \frac{\partial}{\partial x_j}$, $|J|$ denotes the number of elements of $J \subset \{1,\ldots,n\}$, and

(1.2)
$$v(z) \;=\; \frac{\sin\frac{\alpha}{2}(i\beta g + z)}{\sin(\frac{\alpha}{2} z)}.$$

(One has $e^{i\beta\partial_j}\Psi(x_1,\ldots,x_n) = \Psi(x_1,\ldots,x_{j-1}, x_j + i\beta, x_{j+1},\ldots,x_n)$, so $\hat{H}_r$ is a difference operator of order $r$.)

For $g=0$ the coefficients of $\hat{H}_r$ (1.1), (1.2) are trivial ($v=1$) and the difference operators are simultaneously diagonalized by the (Fourier) basis of symmetrized exponentials (or monomial symmetric functions)

(1.3)
$$m_\lambda(x) \;=\; \sum_{\lambda'\in S_n\lambda} e^{i\alpha \sum_{j=1}^n \lambda'_j x_j}, \qquad \lambda \in \Lambda = \{\lambda \in \mathbb{Z}^n \mid \lambda_1 \geq \lambda_2 \geq \cdots \geq \lambda_n\}.$$

For arbitrary $g \geq 0$ a basis of joint eigenfunctions of $\hat{H}_1,\ldots,\hat{H}_n$ follows from work by Macdonald on orthogonal $q$-polynomials in many variables [**12, 13, 14**].

1991 *Mathematics Subject Classification*. Primary 39A70, 47A75, 81Q99; Secondary 33D45.
The author was supported in part by a Shell travel grant.
This paper is in final form and no version of it will be submitted for publication elsewhere.

©0000 American Mathematical Society
0000-0000/00 $1.00 + $.25 per page

One has (assuming $\alpha, \beta > 0$)

$$\hat{H}_r \, \Psi_\lambda \;=\; E_r \left( e^{-\alpha\beta(\lambda_1+\rho_1)}, \ldots, e^{-\alpha\beta(\lambda_n+\rho_n)} \right) \Psi_\lambda, \qquad \lambda \in \Lambda, \tag{1.4}$$

with

$$E_r(t_1, \ldots, t_n) \;=\; \sum_{\substack{J \subset \{1,\ldots,n\} \\ |J|=r}} \prod_{j \in J} t_j, \qquad \rho_j = \tfrac{1}{2} g(n+1-2j), \tag{1.5}$$

and eigenfunctions of the form

$$\Psi_\lambda \;=\; \Delta^{1/2} \, p_\lambda, \qquad \Delta \;=\; \prod_{1 \leq j \neq k \leq n} \frac{(e^{i\alpha(x_j-x_k)}; e^{-\alpha\beta})_\infty}{(e^{-\alpha\beta g} e^{i\alpha(x_j-x_k)}; e^{-\alpha\beta})_\infty}, \tag{1.6}$$

where $p_\lambda$ denotes a trigonometric polynomial that is determined by the conditions

    i. $p_\lambda = m_\lambda + \sum_{\lambda' \in \Lambda, \lambda' < \lambda} c_{\lambda,\lambda'} \, m_{\lambda'}, \qquad c_{\lambda,\lambda'} \in \mathbb{C};$

    ii. $\langle p_\lambda, m_{\lambda'} \rangle_\Delta = 0 \quad \text{if} \quad \lambda' < \lambda.$

Here $\langle \cdot, \cdot \rangle_\Delta$ is the $L^2$ inner product with weight function $\Delta$ (1.6), where, as usual,

$$(a;q)_\infty = \prod_{p=0}^{\infty} (1-aq^p).$$

The (partial) ordering of the cone $\Lambda$ (1.3) is defined by $\lambda' \leq \lambda$ iff $\sum_{j=1}^m \lambda'_j \leq \sum_{j=1}^m \lambda_j$ for $m = 1, \ldots, n$, with equality holding when $m = n$ ($\lambda' < \lambda$ iff $\lambda' \leq \lambda$ and $\lambda' \neq \lambda$). Clearly, $p_\lambda$ amounts to $m_\lambda$ minus its orthogonal projection with respect to $\langle \cdot, \cdot \rangle_\Delta$ on $\text{span}\{m_{\lambda'}\}_{\lambda' \in \Lambda, \lambda' < \lambda}$.

The nonrelativistic limit corresponds to sending the difference step size $\beta$ to zero. Apart from a constant factor the weight function $\Delta$ (1.6) converges to

$$\Delta_0 \;=\; \prod_{1 \leq j < k \leq n} |\sin \tfrac{\alpha}{2}(x_j - x_k)|^{2g}, \tag{1.7}$$

and the eigenfunctions $\Psi_\lambda$ (1.6) go over in $\Psi_{\lambda,0} = \Delta_0^{1/2} p_{\lambda,0}$, where $p_{\lambda,0}$ denotes the trigonometric polynomial determined by conditions i., ii. with the weight function $\Delta$ (1.6) replaced by $\Delta_0$ (1.7). In order to study the transition $\beta \to 0$ for the commuting quantum integrals it is convenient to introduce [20]

$$\tilde{H}_r \;=\; \beta^{-r} \sum_{0 \leq p \leq r} (-1)^p \binom{n-p}{n-r} \hat{H}_p, \qquad r = 1, \ldots, n \tag{1.8}$$

(taking $\hat{H}_0 = 1$ by convention). The operators $\tilde{H}_1, \ldots, \tilde{H}_n$ generate the same algebra of commuting difference operators as $\hat{H}_1, \ldots, \hat{H}_n$, but they have the advantage of having a nontrivial limit when the step size goes to zero. From Eqs. (1.4), (1.8) it follows that

$$\tilde{H}_r \, \Psi_\lambda \;=\; \beta^{-r} \, E_r \left( 1 - e^{-\alpha\beta(\lambda_1+\rho_1)}, \ldots, 1 - e^{-\alpha\beta(\lambda_n+\rho_n)} \right) \Psi_\lambda, \qquad \lambda \in \Lambda, \tag{1.9}$$

which leads for $\beta \to 0$ to

$$\tilde{H}_{r,0} \, \Psi_{\lambda,0} \;=\; \alpha^r \, E_r(\lambda_1+\rho_1, \ldots, \lambda_n+\rho_n) \, \Psi_{\lambda,0}, \qquad \lambda \in \Lambda, \tag{1.10}$$

with

$$\tilde{H}_{r,0} = \lim_{\beta \to 0} \tilde{H}_r = (-i)^r \sum_{\substack{J \subset \{1,\ldots,n\} \\ |J|=r}} \prod_{j \in J} \partial_j + \text{ l.o.} \tag{1.11}$$

(where l.o. stands for terms of lower order in the partials $\partial_j$). Explicit computation of $\tilde{H}_{r,0}$ for $r = 1, 2$ yields

$$\tilde{H}_{1,0} = \frac{1}{i} \sum_{1 \leq j \leq n} \partial_j, \quad \tilde{H}_{2,0} = - \sum_{1 \leq j < k \leq n} \left( \partial_j \partial_k + \tfrac{1}{4} g(g-1) \frac{\alpha^2}{\sin^2 \tfrac{\alpha}{2}(x_j - x_k)} \right). \tag{1.12}$$

The differential operators $\tilde{H}_{1,0}, \ldots, \tilde{H}_{n,0}$ constitute a complete set of commuting quantum integrals for the nonrelativistic Calogero-Sutherland system. Explicit expressions for the integrals can be found in Refs. [3, 21] and [16].

It is well-known that the Calogero-Sutherland model admits a Lie-theoretic generalization in terms of root systems [17]. From this viewpoint the models in this introduction correspond to the root system $A_{n-1}$. In the rest of the paper we are going to focus on the more general case corresponding to the root system $BC_n$. Specifically, we will discuss the diagonalization of a complete set of quantum integrals for a difference version of the quantum Calogero-Sutherland model related to the root system $BC_n$. This difference quantum model may be regarded as a $BC_n$-type generalization of Ruijsenaars' relativistic Calogero-Sutherland system.

*Remarks:* i. The polynomials $p_\lambda$ were introduced by Macdonald [12]. They form a $q$-version of the so-called Jack polynomials [11] (with $q = \exp(-\alpha\beta)$). The transition $p_\lambda \to p_{\lambda,0}$ ($\beta \to 0$) corresponds to the limit $q \to 1$. In this limit Macdonald's polynomials reduce to Jack polynomials.

ii. Conjugation of $\hat{H}_r$ with $\Delta^{1/2}$ eliminates the square roots in the coefficients of the difference operators:

$$\Delta^{-1/2} \hat{H}_r \Delta^{1/2} = \sum_{\substack{J \subset \{1,\ldots,n\} \\ |J|=r}} \prod_{\substack{j \in J \\ k \notin J}} v(x_j - x_k) \, e^{i\beta \sum_{j \in J} \partial_j}. \tag{1.13}$$

The transformed operators (1.13) are diagonalized by the Macdonald polynomials. They are essentially the same as the Macdonald difference operators associated with the fundamental weights of the root system $A_{n-1}$ [13, 14]. (The exact connection is established after restriction to the center of mass plane $x_1 + \cdots + x_n = 0$.)

## 2. Commuting difference operators

The quantum integrals for our model are given by

$$\hat{H}_r = \sum_{\substack{J \subset \{1,\ldots,n\}, |J| \leq r \\ \varepsilon_j = \pm 1, \, j \in J}} U_{J^c,\, r-|J|} \, V_{\varepsilon J;\, J^c}^{1/2} \, e^{i\beta \partial_{\varepsilon J}} \, V_{-\varepsilon J;\, J^c}^{1/2}, \qquad r = 1, \ldots, n, \tag{2.1}$$

with $\partial_{\varepsilon J} = \sum_{j \in J} \varepsilon_j \partial_j$,

$$
(2.2) \quad V_{\varepsilon J; K} = \prod_{j \in J} w(\varepsilon_j x_j) \prod_{\substack{j,j' \in J \\ j<j'}} v(\varepsilon_j x_j + \varepsilon_{j'} x_{j'}) v(\varepsilon_j x_j + \varepsilon_{j'} x_{j'} + i\beta) \\
\times \prod_{\substack{j \in J \\ k \in K}} v(\varepsilon_j x_j + x_k) v(\varepsilon_j x_j - x_k),
$$

(2.3)
$$
U_{K,m} = (-1)^m \sum_{\substack{L \subset K, |L|=m \\ \varepsilon_l = \pm 1, \, l \in L}} \prod_{l \in L} w(\varepsilon_l x_l) \prod_{\substack{l,l' \in L \\ l<l'}} v(\varepsilon_l x_l + \varepsilon_{l'} x_{l'}) v(-\varepsilon_l x_l - \varepsilon_{l'} x_{l'} - i\beta) \\
\times \prod_{\substack{l \in L \\ k \in K \setminus L}} v(\varepsilon_l x_l + x_k) v(\varepsilon_l x_l - x_k),
$$

and

$$
v(z) = \frac{\sin \frac{\alpha}{2}(i\beta g + z)}{\sin(\frac{\alpha}{2} z)},
$$

(2.4)
$$
w(z) = \frac{\sin \frac{\alpha}{2}(i\beta g_0 + z)}{\sin(\frac{\alpha}{2} z)} \frac{\cos \frac{\alpha}{2}(i\beta g_1 + z)}{\cos(\frac{\alpha}{2} z)} \frac{\sin \frac{\alpha}{2}(i\beta g_0' + \frac{i\beta}{2} + z)}{\sin \frac{\alpha}{2}(\frac{i\beta}{2} + z)} \frac{\cos \frac{\alpha}{2}(i\beta g_1' + \frac{i\beta}{2} + z)}{\cos \frac{\alpha}{2}(\frac{i\beta}{2} + z)}.
$$

It is clear that, just as the quantum integrals for the models considered in the introduction, these difference operators are invariant with respect to permutations of the particle positions $x_1, \ldots, x_n$. However, the operators are no longer translational invariant; instead, we now have invariance with respect to the reflections $x_j \to -x_j$, $j = 1, \ldots, n$. Furthermore, apart from the potential $v$, which governs the interaction between the particles, there is now a second potential $w$ modeling an external field. From now on we will always assume that

$$
(2.5) \qquad \alpha, \beta > 0 \quad \text{and} \quad g, g_0^{(\prime)}, g_1^{(\prime)} \geq 0.
$$

These constraints ensure that our quantum integrals $\hat{H}_r$ (2.1)-(2.4) are self-adjoint. (One has $\overline{v(x_j)} = v(-x_j)$, $\overline{w(x_j)} = w(-x_j)$, so $\overline{V}_{\varepsilon J; K} = V_{-\varepsilon J; K}$ and $\overline{U}_{K,m} = U_{K,m}$.)

The (weight) function that transforms away the square roots in the coefficients of the difference operators (cf. Remark *ii.* in the introduction) now reads

$$
(2.6) \qquad \Delta = \prod_{\substack{1 \leq j < j' \leq n \\ \varepsilon, \varepsilon' = \pm 1}} d_v(\varepsilon x_j + \varepsilon' x_{j'}) \prod_{\substack{1 \leq j \leq n \\ \varepsilon = \pm 1}} d_w(\varepsilon x_j),
$$

with

$$
d_v(z) = \frac{(e^{i\alpha z}; e^{-\alpha\beta})_\infty}{(e^{-\alpha\beta g} e^{i\alpha z}; e^{-\alpha\beta})_\infty},
$$

$$d_w(z) = \frac{(e^{i\alpha z}; e^{-\alpha\beta})_\infty}{(e^{-\alpha\beta g_0} e^{i\alpha z}; e^{-\alpha\beta})_\infty} \frac{(-e^{i\alpha z}; e^{-\alpha\beta})_\infty}{(-e^{-\alpha\beta g_1} e^{i\alpha z}; e^{-\alpha\beta})_\infty}$$
$$\times \frac{(e^{i\alpha(i\beta/2+z)}; e^{-\alpha\beta})_\infty}{(e^{-\alpha\beta g'_0} e^{i\alpha(i\beta/2+z)}; e^{-\alpha\beta})_\infty} \frac{(-e^{i\alpha(i\beta/2+z)}; e^{-\alpha\beta})_\infty}{(-e^{-\alpha\beta g'_1} e^{i\alpha(i\beta/2+z)}; e^{-\alpha\beta})_\infty}.$$

Conjugating with $\Delta^{1/2}$ results in

$$(2.7) \quad \hat{\mathcal{H}}_r = \Delta^{-1/2} \hat{H}_r \Delta^{1/2} = \sum_{\substack{J \subset \{1,\ldots,n\}, |J| \leq r \\ \varepsilon_j = \pm 1, j \in J}} U_{J^c, r-|J|} V_{\varepsilon J; J^c} e^{i\beta \partial_{\varepsilon J}}.$$

The commutativity of $\hat{H}_1, \ldots, \hat{H}_n$ (or equivalently $\hat{\mathcal{H}}_1, \ldots, \hat{\mathcal{H}}_n$) is by no means immediate from the above expressions, but it follows from the fact that (as we will see below) the operators have a basis of joint eigenfunctions.

## 3. Triangularity

The first step in the diagonalization process consists of demonstrating that the transformed difference operators $\hat{\mathcal{H}}_r$ (2.7), (2.2)-(2.4) map the space of permutation invariant, even, trigonometric polynomials into itself. This space is spanned by the basis of $W$-invariant monomials

(3.1)
$$m_\lambda(x) = \sum_{\lambda' \in W\lambda} e^{i\alpha \sum_{j=1}^n \lambda'_j x_j}, \quad \lambda \in \Lambda = \{\lambda \in \mathbb{Z}^n \mid \lambda_1 \geq \lambda_2 \geq \cdots \geq \lambda_n \geq 0\},$$

where $W$ is the group generated by permutations and sign flips ($W \cong S_n \ltimes (\mathbb{Z}_2)^n$). When acting with $\hat{\mathcal{H}}_r$ on a monomial $m_\lambda$ we obtain a $W$-invariant trigonometric function that is rational in $\exp(\pm i\alpha x_j)$, $j = 1, \ldots, n$:

$$(3.2) \quad \hat{\mathcal{H}}_r m_\lambda = \sum_{\substack{J_+ \cap J_- = \emptyset, |J_+ \cup J_-| = r \\ \varepsilon_j = \pm 1, j \in J_+ \cup J_-}} V_{\varepsilon J_+, \varepsilon J_-; J_+^c \cap J_-^c} \, m_\lambda(x + i\beta \sum_{j \in J_+} \varepsilon_j e_j),$$

with $e_j$ denoting the $j$th unit vector of the standard basis of $\mathbb{R}^n$ and

(3.3)
$$V_{\varepsilon J_+, \varepsilon J_-; K} = (-1)^{|J_-|} \prod_{j \in J_+ \cup J_-} w(\varepsilon_j x_j) \prod_{\substack{j \in J_+ \cup J_- \\ k \in K}} v(\varepsilon_j x_j + x_k) v(\varepsilon_j x_j - x_k)$$
$$\times \prod_{\substack{j,j' \in J_+ \\ j < j'}} v(\varepsilon_j x_j + \varepsilon_{j'} x_{j'}) v(\varepsilon_j x_j + \varepsilon_{j'} x_{j'} + i\beta)$$
$$\times \prod_{\substack{j,j' \in J_- \\ j < j'}} v(\varepsilon_j x_j + \varepsilon_{j'} x_{j'}) v(-\varepsilon_j x_j - \varepsilon_{j'} x_{j'} - i\beta)$$
$$\times \prod_{\substack{j \in J_+ \\ j' \in J_-}} v(\varepsilon_j x_j + x_{j'}) v(\varepsilon_j x_j - x_{j'}).$$

To prove that $\hat{\mathcal{H}}_r m_\lambda$ is actually polynomial in $\exp(\pm i\alpha x_1), \ldots, \exp(\pm i\alpha x_n)$, it suffices to demonstrate that, viewed as a function of $x_j$ (with all other variables fixed in general position), the r.h.s. of (3.2) is regular. Thus, we must show that the (generically) simple poles in the terms of $\hat{\mathcal{H}}_r \, m_\lambda$, which are caused by the zeros in the denominators of the potentials $v, w$ (2.4), all cancel each other in (3.2). The poles are congruent (mod $2\pi/\alpha$) to one of

$$
\begin{aligned}
x_j &\sim 0 \quad (+\tfrac{\pi}{\alpha}), & &\text{type I}, \\
x_j &\sim \pm x_{j'}, & j' &\neq j, \quad \text{type II}, \\
x_j &\sim \pm i\beta/2 \quad (+\tfrac{\pi}{\alpha}), & &\text{type III}, \\
x_j &\sim \pm i\beta \pm x_{j'}, & j' &\neq j, \quad \text{type IV}.
\end{aligned}
$$

For the poles of type I and type III it is enough to consider only those congruent to $0$ and $\pm i\beta/2$, respectively. This is because a translation over the half-period $\pi/\alpha$ is equivalent to an interchange of parameters $g_0^{(\prime)} \leftrightarrow g_1^{(\prime)}$. Let us first consider the types I and II. The residue at $x_j = 0$ vanishes because (3.2) is even in $x_j$. The vanishing of the residues at the poles of type II is also immediate from the $W$-invariance: we first use the permutation symmetry to conclude the vanishing of the residue at $x_j = x_{j'}$, and then the invariance under sign flips to conclude the same for the residue at $x_j = -x_{j'}$. For the remaining two cases (type III and type IV) we may assume that all signs are negative. This corresponds to considering terms in (3.2) with $\varepsilon_j, \varepsilon_{j'} = +1$. The poles of type III arise in the terms with $j$ either in $J_+$ or in $J_-$. By combining each term with the term corresponding to the index sets obtained by bringing $j$ from $J_+$ to $J_-$ or vice versa, it is seen that the residues cancel each other pairwise. The situation for poles of type IV is very similar; these poles only occur in terms with both $j, j'$ either in $J_+$ or in $J_-$. By bringing $j, j'$ from $J_+$ to $J_-$ or vice versa, one infers that the residues of the corresponding terms also cancel in pairs.

Now that it has been verified that the r.h.s. of Eq. (3.2) is regular in $x_j$, $j = 1, \ldots, n$, we may conclude that $\hat{\mathcal{H}}_r \, m_\lambda$ is a $W$-invariant trigonometric polynomial. Hence, it can be written as a finite linear combination of the monomials (3.1). It turns out that the expansion of $\hat{\mathcal{H}}_r \, m_\lambda$ in monomials is of a very specific form:

$$(3.4) \qquad \hat{\mathcal{H}}_r \, m_\lambda \;=\; \sum_{\substack{\lambda' \in \Lambda \\ \lambda' \leq \lambda}} [\hat{\mathcal{H}}_r]_{\lambda, \lambda'} \; m_{\lambda'}, \qquad \text{with} \quad [\hat{\mathcal{H}}_r]_{\lambda, \lambda'} \in \mathbb{C},$$

where the (partial) ordering of the cone $\Lambda$ (3.1) labeling the monomial basis is defined by

$$(3.5) \qquad \lambda' \leq \lambda \qquad \text{iff} \qquad \sum_{1 \leq j \leq m} \lambda'_j \leq \sum_{1 \leq j \leq m} \lambda_j \quad \text{for} \quad m = 1, \ldots, n.$$

In other words: the operator $\hat{\mathcal{H}}_r$ is triangular with respect to the partially ordered basis of $W$-invariant monomials $\{m_\lambda\}_{\lambda \in \Lambda}$.

The triangularity of $\hat{\mathcal{H}}_r$ hinges on the fact that the potentials $v(z), w(z)$ (2.4) have constant asymptotics for $\text{Im}(z) \to \pm\infty$:

$$(3.6) \qquad \lim_{\text{Im}(z) \to \pm\infty} v(z) \;=\; e^{\pm \alpha \beta g/2}, \qquad \lim_{\text{Im}(z) \to \pm\infty} w(z) \;=\; e^{\pm \alpha \beta (g_0 + g_1 + g'_0 + g'_1)/2}.$$

To see this we put $x = iRy$, with $y \in \mathbb{R}^n$ and $y_1 > y_2 > \cdots > y_n > 0$, and combine the asymptotics

$$(3.7) \qquad m_{\lambda'}(iRy) \;=\; e^{\alpha R \sum_{j=1}^n \lambda'_j y_j} \, (1 + o(1)) \qquad \text{for} \quad R \to \infty$$

with Eq. (3.6) to conclude that

$$(3.8) \quad (\hat{\mathcal{H}}_r m_\lambda)(iRy) = O(e^{\alpha R \sum_{j=1}^n \lambda_j y_j}) \quad \text{for} \quad R \to \infty.$$

By comparing the asymptotics (3.7) and (3.8), and using the fact that

$$(3.9) \quad \lambda' \leq \lambda \quad \text{iff} \quad \sum_{1 \leq j \leq n} \lambda'_j y_j \leq \sum_{1 \leq j \leq n} \lambda_j y_j \quad \forall y \in \mathbb{R}^n \text{ with } y_1 > \cdots > y_n > 0,$$

one infers that in the expansion of $\hat{\mathcal{H}}_r m_\lambda$ in monomials $m_{\lambda'}$ only the terms with $\lambda' \leq \lambda$ may occur with a nonzero coefficient.

## 4. Eigenfunctions and eigenvalues

Let $p_\lambda$ be $m_\lambda$ minus its orthogonal projection onto $\text{span}\{m_{\lambda'}\}_{\lambda' \in \Lambda, \lambda' < \lambda}$ with respect to the inner product $\langle \cdot, \cdot \rangle_\Delta$ (i.e., the $L^2$ inner product with weight function $\Delta$ (2.6)). Clearly $p_\lambda$ is uniquely determined as the monic polynomial in $\text{span}\{m_{\lambda'}\}_{\lambda' \in \Lambda, \lambda' \leq \lambda}$ that is orthogonal to $\text{span}\{m_{\lambda'}\}_{\lambda' \in \Lambda, \lambda' < \lambda}$. The triangularity of $\hat{\mathcal{H}}_r$ ensures that $\hat{\mathcal{H}}_r p_\lambda$ lies also in $\text{span}\{m_{\lambda'}\}_{\lambda' \in \Lambda, \lambda' \leq \lambda}$, and the symmetry of $\hat{\mathcal{H}}_r$ with respect to $\langle \cdot, \cdot \rangle_\Delta$ implies that $\hat{\mathcal{H}}_r p_\lambda$ is orthogonal to $\text{span}\{m_{\lambda'}\}_{\lambda' \in \Lambda, \lambda' < \lambda}$. (One has $\langle \hat{\mathcal{H}}_r p_\lambda, m_{\lambda'} \rangle_\Delta = \langle p_\lambda, \hat{\mathcal{H}}_r m_{\lambda'} \rangle_\Delta = 0$ if $\lambda' < \lambda$.) Hence, $\hat{\mathcal{H}}_r p_\lambda$ must be proportional to $p_\lambda$, i.e., $p_\lambda$ is an eigenfunction of $\hat{\mathcal{H}}_r$.

Because the basis transformation $\{m_\lambda\}_{\lambda \in \Lambda} \to \{p_\lambda\}_{\lambda \in \Lambda}$ does not affect the matrix elements on the diagonal, the corresponding eigenvalues are given by (using Eqs. (3.4), (3.7) and (3.9), together with the fact that in the last equation the inequalities are strict if $\lambda' < \lambda$):

$$(4.1) \quad [\hat{\mathcal{H}}_r]_{\lambda,\lambda} = \lim_{R \to \infty} e^{-\alpha R \sum_{j=1}^n \lambda_j y_j} (\hat{\mathcal{H}}_r m_\lambda)(iRy).$$

To evaluate this limit we use Eqs. (3.2), (3.3) and (3.6), (3.7) to obtain

$$[\hat{\mathcal{H}}_r]_{\lambda,\lambda} = \sum_{\substack{J_+ \cap J_- = \emptyset, |J_+ \cup J_-|=r \\ \varepsilon_j = \pm 1, j \in J_+ \cup J_-}} (-1)^{|J_-|} \prod_{j \in J_+} e^{\alpha\beta\varepsilon_j(\lambda_j + \rho_j)} \prod_{j \in J_-} e^{\alpha\beta\varepsilon_j(g_0 + g_1 + g'_0 + g'_1)/2}$$

$$\times \prod_{\substack{j \in J_-, k \in J_+^c \cap J_-^c \\ j < k}} e^{\alpha\beta\varepsilon_j g}$$

$$= 2^r \sum_{\substack{J \subset \{1,\ldots,n\} \\ 0 \leq |J| \leq r}} (-1)^{r-|J|} \prod_{j \in J} \text{ch}\,\alpha\beta(\lambda_j + \rho_j) \sum_{r \leq l_1 \leq \cdots \leq l_{r-|J|} \leq n} \text{ch}\,\alpha\beta\rho_{l_1} \cdots \text{ch}\,\alpha\beta\rho_{l_{r-|J|}},$$

with $\rho_j = (n-j)\,g + (g_0 + g_1 + g'_0 + g'_1)/2$. Transforming back to Lebesgue measure leads to the eigenfunctions of $\hat{H}_1, \ldots, \hat{H}_n$ (2.1)-(2.4).

Let us summarize the final result:

$$\hat{H}_r \Psi_\lambda = E_r\left(\text{ch}\,\alpha\beta(\lambda_1 + \rho_1), \ldots, \text{ch}\,\alpha\beta(\lambda_n + \rho_n); \text{ch}(\alpha\beta\rho_r), \ldots, \text{ch}(\alpha\beta\rho_n)\right) \Psi_\lambda,$$

with

$$E_r(t_1,\ldots,t_n;p_r,\ldots,p_n) = 2^r \sum_{\substack{J\subset\{1,\ldots,n\}\\ 0\leq |J|\leq r}} (-1)^{r-|J|} \prod_{j\in J} t_j \sum_{r\leq l_1\leq\cdots\leq l_{r-|J|}\leq n} p_{l_1}\cdots p_{l_{r-|J|}}$$

and eigenfunctions of the form

$$\Psi_\lambda = \Delta^{1/2} p_\lambda, \qquad \lambda \in \Lambda = \{\lambda\in\mathbb{Z}^n \mid \lambda_1\geq \lambda_2\geq\cdots\geq\lambda_n\geq 0\},$$

where $\Delta$ is taken from (2.6) and $p_\lambda$ denotes the trigonometric polynomial determined by the conditions

i. $p_\lambda = m_\lambda + \sum_{\lambda'\in\Lambda,\lambda'<\lambda} c_{\lambda,\lambda'} m_{\lambda'}, \quad c_{\lambda,\lambda'}\in\mathbb{C};$

ii. $\langle p_\lambda, m_{\lambda'}\rangle_\Delta = 0 \quad\text{if}\quad \lambda'<\lambda$

with the partial ordering of the cone $\Lambda$ as defined in (3.5).

*Remark:* For $r=1$ the operator $\hat{\mathcal{H}}_r$ (2.7), (2.2)-(2.4) becomes

$$(4.2)\qquad \hat{\mathcal{H}}_1 = \sum_{\substack{1\leq j\leq n\\ \varepsilon=\pm 1}} w(\varepsilon\, x_j) \prod_{k\neq j} v(\varepsilon\, x_j + x_k)\, v(\varepsilon\, x_j - x_k)\, (e^{i\beta\varepsilon\partial_j} - 1).$$

This operator and its eigenfunctions $p_\lambda$, $\lambda\in\Lambda$, were introduced by Koornwinder [10] as a generalization of Macdonald's construction [13] for the root system $BC_n$. All Macdonald polynomials associated with (admissible pairs of) the classical root systems (i.e., $A_{n-1}$, $B_n$, $C_n$, $D_n$, and $BC_n$) can be seen as special cases of the Koornwinder polynomials $\{p_\lambda\}_{\lambda\in\Lambda}$ [6]. In the case of one variable ($n=1$) Koornwinder's polynomials reduce to Askey-Wilson polynomials [1].

## 5. Concluding remarks

**5.1. Trigonometric identities.** In Ref. [6] we studied the diagonalization of difference operators $\hat{\mathcal{H}}_r$ (2.7) with $V_{\varepsilon J;K}$ as in (2.2) and $U_{K,m}$ given by

$$(5.1)\quad U_{K,m} = \sum_{\substack{\emptyset\subsetneq L_1\subsetneq\cdots\subsetneq L_p\subset K,\, 1\leq p\leq m\\ |L_p|=m,\ \varepsilon_l=\pm 1,\ l\in L_p}} (-1)^p\ V_{\varepsilon L_1;K\setminus L_1} V_{\varepsilon(L_2\setminus L_1);K\setminus L_2}\cdots V_{\varepsilon(L_p\setminus L_{p-1});K\setminus L_p}$$

(where $v$ and $w$ are again taken from (2.4)). We found that these difference operators are simultaneously diagonalized by Koornwinder's multivariable Askey-Wilson polynomials, with the same eigenvalues as the operators in the present paper (i.e., with $U_{K,m}$ given by (2.3)). Because the difference operators are uniquely determined by their eigenvalues [6], it follows that the r.h.s. of Eq. (2.3) equals the r.h.s. of Eq. (5.1). This equality hinges on nontrivial functional identities for the trigonometric function $v$ (2.4) [7].

**5.2. The limit $\beta\to 0$.** For $\beta\to 0$, one has [6]

$$\hat{H}_r = \beta^{2r}\, \hat{H}_{r,0} + o(\beta^{2r}), \qquad \hat{H}_{r,0} = (-1)^r \sum_{\substack{J\subset\{1,\ldots,n\}\\ |J|=r}} \prod_{j\in J} \partial_j^2 + \text{l.o.},$$

with

$$\hat{H}_{r,0}\,(\Delta_0^{\frac{1}{2}}\,p_{\lambda,0}) =$$
$$\alpha^{2r} E_r((\lambda_1+\rho_1)^2,\ldots,(\lambda_n+\rho_n)^2;\rho_r^2,\ldots,\rho_n^2)\,\Delta_0^{\frac{1}{2}}\,p_{\lambda,0}, \qquad \lambda \in \Lambda,$$

where $p_{\lambda,0}$ denotes a multivariable generalization of the Jacobi polynomials [4, 9] determined by the conditions i., ii. with weight function

$$\Delta_0 = \prod_{1\leq j<k\leq n} |\sin\tfrac{\alpha}{2}(x_j+x_k)\,\sin\tfrac{\alpha}{2}(x_j-x_k)|^{2g}\prod_{1\leq j\leq n}|\sin(\tfrac{\alpha}{2}x_j)|^{2\tilde{g}_0}|\cos(\tfrac{\alpha}{2}x_j)|^{2\tilde{g}_1},$$

and $\tilde{g}_0 = g_0 + g_0'$, $\tilde{g}_1 = g_1 + g_1'$. The differential operators $\hat{H}_{1,0},\ldots,\hat{H}_{n,0}$ constitute a complete set of quantum integrals for the Calogero-Sutherland system related to the root system $BC_n$. For $r=1$ the operator reads

$$\hat{H}_{1,0} = -\sum_{1\leq j\leq n}\partial_j^2 + \tfrac{1}{2}g(g-1)\alpha^2\sum_{1\leq j<k\leq n}\left(\frac{1}{\sin^2\tfrac{\alpha}{2}(x_j+x_k)} + \frac{1}{\sin^2\tfrac{\alpha}{2}(x_j-x_k)}\right)$$
$$+\tfrac{1}{4}\alpha^2\sum_{1\leq j\leq n}\left(\frac{\tilde{g}_0(\tilde{g}_0-1)}{\sin^2(\tfrac{\alpha}{2}x_j)}+\frac{\tilde{g}_1(\tilde{g}_1-1)}{\cos^2(\tfrac{\alpha}{2}x_j)}\right) \;+\text{constant}.$$

Explicit expressions for the higher-order quantum integrals can be found in Refs. [4, 16].

### 5.3. Differential- and difference-reflection operators.
For arbitrary root system the *existence* of a complete set of quantum integrals for the Calogero-Sutherland system can be proved using a construction involving certain differential-reflection operators known as 'Dunkl operators' [9]. This approach does not yield explicit expressions for the quantum integrals, but it does provide so-called shift operators. With these shift operators well-known conjectures due to Macdonald for the values of normalization constants $\langle p_{\lambda,0}, p_{\lambda,0}\rangle_{\Delta_0}$ can be proved [18]. A partial generalization of this method to the Macdonald polynomials related to reduced root systems was presented recently by Cherednik [2]. Specifically, it follows from Cherednik's work that for special values of the parameters $g_0^{(\prime)}$, $g_1^{(\prime)}$ the *existence* of a complete set of quantum integrals for $\hat{\mathcal{H}}_1$ (4.2) can also be shown with aid of difference-reflection operators, and that, in these special cases, Macdonald's $q$-version of the normalization conjectures for $\langle p_\lambda, p_\lambda\rangle_\Delta$ can be proved.

### 5.4. Quantum group interpretations.
It is expected that for certain special values of the parameters the Koornwinder-Askey-Wilson polynomials can be seen as zonal spherical functions on quantum versions of homogeneous spaces of classical simple Lie groups. Our commuting difference operators should then correspond to radial reductions of the Casimir elements. For the root system $A_{n-1}$ such a quantum group interpretation of the Macdonald polynomials was given by Noumi for $g=1/2$ and $g=2$ [15]. A different approach, valid for arbitrary values of $g$, was developed in Ref. [8].

### 5.5. Elliptic generalizations.
It is well-known that the trigonometric Calogero-Sutherland system has an integrable generalization involving elliptic potentials [17]. Explicit quantum integrals for the elliptic version can be found in [16] (differential case, all classical root systems) and [19] (difference case, root system $A_{n-1}$). A partial generalization to elliptic potentials of $\hat{H}_r$ (2.1)-(2.4) was presented in [5].

DEPARTMENT OF MATHEMATICS AND COMPUTER SCIENCE, UNIVERSITY OF AMSTERDAM, PLANTAGE MUIDERGRACHT 24, 1018 TV AMSTERDAM, THE NETHERLANDS